\documentclass[prl,aps,twocolumn,10pt,nofootinbib]{revtex4-2}
\usepackage{amsmath}
\usepackage{hyperref}
\usepackage{orcidlink}
\usepackage[dvipsnames]{xcolor} 
\usepackage[english]{babel} 
\usepackage[autostyle, english = american]{csquotes} 
  \MakeOuterQuote{"} 
\usepackage{graphicx}
  \graphicspath{{figures}}

\newcommand{\musig}{\mu^{\textrm{(signal)}}}
\newcommand{\Ncrit}{N_{\textrm{crit}}}

\usepackage{pdfpages}
\usepackage{etoolbox} 

\makeatletter
\patchcmd{\@outputpage@head}{\@ifx{\LS@rot\@undefined}{}{\LS@rot}}{}{}{}
\makeatother

\begin{document}

\title[AIC issues]{Information Criteria Fail for Dynamical Systems: Sampling Rate and Dimension Dependence}

\author{Kumar Utkarsh}
\email{krutkarsh@u.northwestern.edu.}

\author{Daniel M.~Abrams, \orcidlink{0000-0002-6015-8358}}%
\email{dmabrams@northwestern.edu}
\affiliation{%
Department of Engineering Sciences and Applied Mathematics, McCormick School of Engineering and Applied Science, Northwestern University, 2145 Sheridan Road, Evanston, Illinois 60208, USA
}%

\begin{abstract}
    Information criteria such as Akaike's (AIC) and Bayes' (BIC) are widely used for model selection in physics and beyond, quantifying the tradeoff between model complexity and goodness-of-fit to enforce parsimony. However, their derivation assumes uncorrelated samples, an assumption systematically violated by dynamical systems data. Here, through analysis of simple but representative dynamical models---exponential decay, harmonic oscillation, and chaos---we demonstrate that model selection depends sensitively on sampling rate and system dimensionality. We derive explicit formulas predicting when standard information criteria fail that should be adaptable to many real-world scenarios, enabling experimentalists to design sampling protocols that avoid pathological regimes.  
\end{abstract}

\maketitle

\section{Introduction}
\label{sec:intro}

Information criteria, especially the Akaike Information Criterion (AIC) \cite{akaike1974new} and Bayesian Information Criterion (BIC) \cite{schwarz1978estimating}, are fundamental tools for model selection in many fields including physics, statistics, and machine learning \cite{burnham2002model, konishi2008, vrieze2012model, claeskens2008model}. Their theoretical foundation lies in approximating the Kullback-Leibler (KL) divergence \cite{kullback1951information} between the true data-generating process and the candidate models, and, in a deeper sense, in information theory's connection to mathematical models \cite{transtrum2010nonlinear, transtrum2015perspective, gutenkunst2007universally}. 

These criteria have found extensive application across physics (see, e.g., \cite{usami2003accuracy, neil2024improved, anisimo2011spectro, hippel2025qcd, figueroa2024cosmo}) and, indeed across all of quantitative science (see, e.g., \cite{hamilton1994time, posada2008jmodeltest, aho2014model}. However, their application to dynamical systems with temporal correlations has received limited theoretical scrutiny. While it is well established that information criteria assume independent observations \cite{claeskens2008model}---an assumption usually violated by dynamical systems data---a tractable predictive framework for experimental design remains absent. Existing approaches either restrict themselves to a particular class of models \cite{jones2011bayesian, friston2003dynamic} or require computationally intensive simulation-based methods \cite{shang2008bootstrap,toni2009approximate}.
What is missing are explicit scaling laws that enable experimentalists to understand, \textit{a priori}, which sampling protocols will yield reliable model selection and which will lead to systematic failures.

Here, we present a series of simple examples to illustrate two key problems:\\
 \textbf{Sampling rate dependence.} The model selected by an information criterion can depend sensitively on the sampling rate of the observations. This is counterintuitive, as one might expect more data (higher sampling rates) to always improve model selection. However, both under-sampling and over-sampling can lead to incorrect conclusions because information criteria assume independent observations \cite{claeskens2008model}. This problem is particularly acute in physics, where experimenters often have significant control over sampling parameters but may unknowingly choose rates that lead to systematic model misidentification \cite{white2016limitations, ding2006granger}.\\
 \textbf{Dimension dependence.} The reliability of model selection depends on the system dimension. Surprisingly, this represents a ``dimensionality blessing'' rather than the usual curse---higher-dimensional systems with many independent realizations can overcome complexity penalties through collective evidence \cite{schwarz1978estimating,stone1977asymptotic}. However, this blessing depends critically on how the number of parameters scales with dimension, an issue that has received limited attention in the context of collective dynamical systems \cite{transtrum2015perspective}.

These pathologies are particularly concerning because they can lead to systematic misidentification of physical mechanisms. For instance, genuine relaxation dynamics might be dismissed as noise due to poor sampling choices, or tractable models for behavior might be overlooked in low-dimensional observations \cite{daniels2015automated}. The implications extend beyond model selection to parameter estimation and uncertainty quantification \cite{machta2013parameter}. 

Our analysis provides complete characterization for first-order linear systems and harmonic oscillators---ubiquitous dynamical motifs appearing throughout physics. The resulting scaling laws offer immediate practical guidance: experimentalists can evaluate whether their sampling protocols fall in reliable or pathological regimes, and adjust accordingly before data collection. Unlike previous studies that focus on comparing different information criteria (e.g., \cite{aho2014model}), we demonstrate that the problems are intrinsic to the application of likelihood-based methods to temporally correlated data, regardless of which specific criterion is used\footnote{We suspect that similar problems also exist in practice for non-likelihood methods like minimum description length (MDL).}.  Our goal is to establish a tractable predictive framework for experimental design to mitigate these problems.

\section{Sampling Rate Dependence}
\label{sec:sampling}
\noindent \textbf{Equilibrating Systems:}
Consider a simple deterministic decay model governed by
\begin{equation}
    \dot{x} = -\lambda x,
    \label{eq:decay}
\end{equation}
where $\lambda > 0$ is the decay rate and the overdot indicates a time derivative.  Assume an initial condition $x(0)=x_0$ and normally distributed observation noise $\eta \sim \mathcal{N}(\mu,\sigma^2)$, where $\mu$ and $\sigma^2$ are the noise mean and variance, respectively.  Then the observations should satisfy $x(t) \sim x_0 e^{-\lambda t} +\eta$.

\smallskip
\noindent{\textit{Fixed total data:}} Suppose we sample $x(t)$ at a fixed sampling frequency $f$ (i.e., at discrete times with interval $\Delta t = 1/f$), collecting $M$ total data points.  We wish to answer the key question: \textit{are these data more consistent with the decay model in Eq.~\eqref{eq:decay}, or are they more consistent with a null model that assumes pure noise?  }

In the low sampling frequency limit $f \to 0^+$ (large $\Delta t$), the sampled data resemble independent Gaussian draws, leading to selection of the ``pure noise'' model.  As sampling frequency grows, however, temporal correlations emerge and the decay model is eventually favored.  In the high sampling frequency limit $f \to \infty$ (small $\Delta t$), however, the sampled data again resemble independent Gaussian draws, and once again the pure noise (null) model is selected.

The crossovers between these different selection regimes depend on the decay rate $\lambda$, initial condition $x_0$, noise level $\sigma$, and number of samples $M$. For the low-frequency crossover, sparse sampling fails to capture sufficient decay structure when $f \ll \lambda$. In this regime, only the first $\sim f/\lambda$ data points contain significant signal, and the signal variance relative to noise variance becomes too small to justify the extra complexity penalty. Detailed analysis (see Appendix~\ref{app:freqcrossover}) yields the lower crossover frequency:
\begin{equation}
f_c^{(1)} = \frac{8\lambda\sigma^2}{x_0^2}.
\label{eq:decaycrosslow}
\end{equation}

For the high-frequency crossover, the observation window $t_{\max} = M/f$ becomes too short to capture decay dynamics when $f \gg \lambda$. In this regime, the signal appears nearly constant over the sampling period, with variance scaling as $x_0^2\Lambda^2/12$ where $\Lambda = \lambda M/f \ll 1$ is the number of decay timescales observed. Balancing this against the AIC complexity penalty yields the upper crossover frequency:
\begin{equation}
f_c^{(2)} = \frac{M^{3/2}\lambda x_0^2}{4\sqrt{3}\sigma}.
\label{eq:decaycrosshigh}
\end{equation}

\begin{figure}[t]
    \centering
        \includegraphics[width=\linewidth]{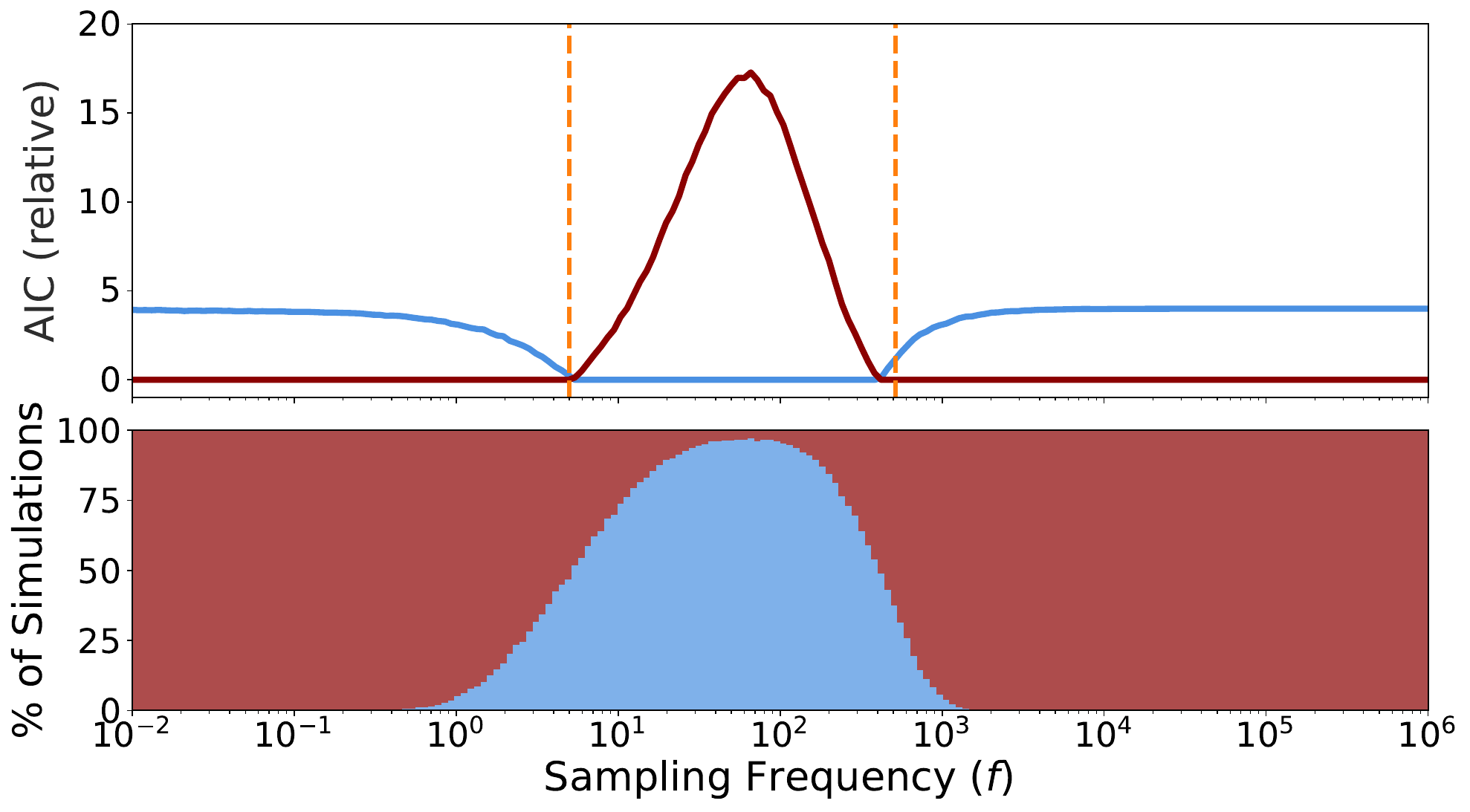}
    \caption{\textbf{Equilibration: Selected model varies with sampling rate.} Top: Comparison of AIC values for true deterministic decay model (blue) and pure noise null model (red) as a function of sampling frequency $f$. Theoretical crossover points are indicated by vertical orange dashed lines: $f_c^{(1)} = 8\lambda\sigma^2/x_0^2 = 5.0$ and $f_c^{(2)} = M^{3/2}\lambda x_0^2/(4\sqrt{3}\sigma) \approx 2000$.  Bottom: Proportion of trials in which each model is selected as a function of sampling frequency $f$ (1000 trials). Note that the noise model is selected at high and low frequencies (left and right on graph). Parameters: $(x_0, \lambda, \mu, \sigma, M) = (1, 0.1, 0, 2.5, 2000)$. Points in upper panel represent means over trials at each frequency. 
    }
    \label{fig:decaymodelcomp}
\end{figure}

So for $f_c^{(1)} \lesssim f \lesssim f_c^{(2)}$ we expect the ``true'' decay model to be selected, but for sampling frequencies much higher or lower we expect the noise model to be selected.  Fig.~\ref{fig:decaymodelcomp} demonstrates numerically that this is indeed the result when using AIC.

This analysis clearly demonstrates the counter-intuitive result that both under-sampling and over-sampling can lead to incorrect model selection, even when the true model is included in the candidate set. This suggests that optimal experimental design for dynamical systems must carefully consider intrinsic time scales of the system being studied.

\noindent \textit{Fixed sampling frequency:} 
Repeating the above numerical experiment but now fixing sampling frequency $f$ instead of the total number of data points collected $M$, we find a single crossover. This crossover occurs when the observation window $t_{\max} = M/f$ becomes long enough to capture sufficient decay structure. From the high-frequency crossover analysis, this happens when $f = M^{3/2}\lambda x_0/(4\sqrt{3}\sigma)$, which gives a critical sample size of
\begin{equation}
    M_c = \left(\frac{4\sqrt{3}\sigma f}{\lambda x_0^2}\right)^{2/3}.
    \label{eq:mcrossoverdecay}
\end{equation}
Note that the above crossover applies at high sampling frequencies; deviations occur as sampling frequency decreases. Both fixed $M$ and fixed $f$ cases reveal a central principle: \textit{recovery of the true model occurs only when the sampling resolution and the total observation window are jointly matched to the intrinsic time scale of the system}.

\smallskip

\noindent \textbf{Oscillatory systems:} Beyond equilibration, perhaps the next most common qualitative dynamical behavior is oscillation. Here we consider as a prototypical example the simple harmonic oscillator (SHO); our arguments should apply well to any stable or neutrally stable limit cycle in some capacity. Consider the ordinary differential equation (ODE)
\begin{equation}
    \ddot{x}(t) = -\omega^2 x(t).
\end{equation}
Oscillatory solutions are parametrized by amplitude $A$ and initial phase $\phi$ as $x(t) = A \cos(\omega t + \phi)$. Without loss of generality, we consider the case $\phi=0$. Sampling $x(t)$ at different rates alters the statistical properties of the resulting time series.

\begin{figure}[t]
    \centering
    \includegraphics[width=\linewidth]{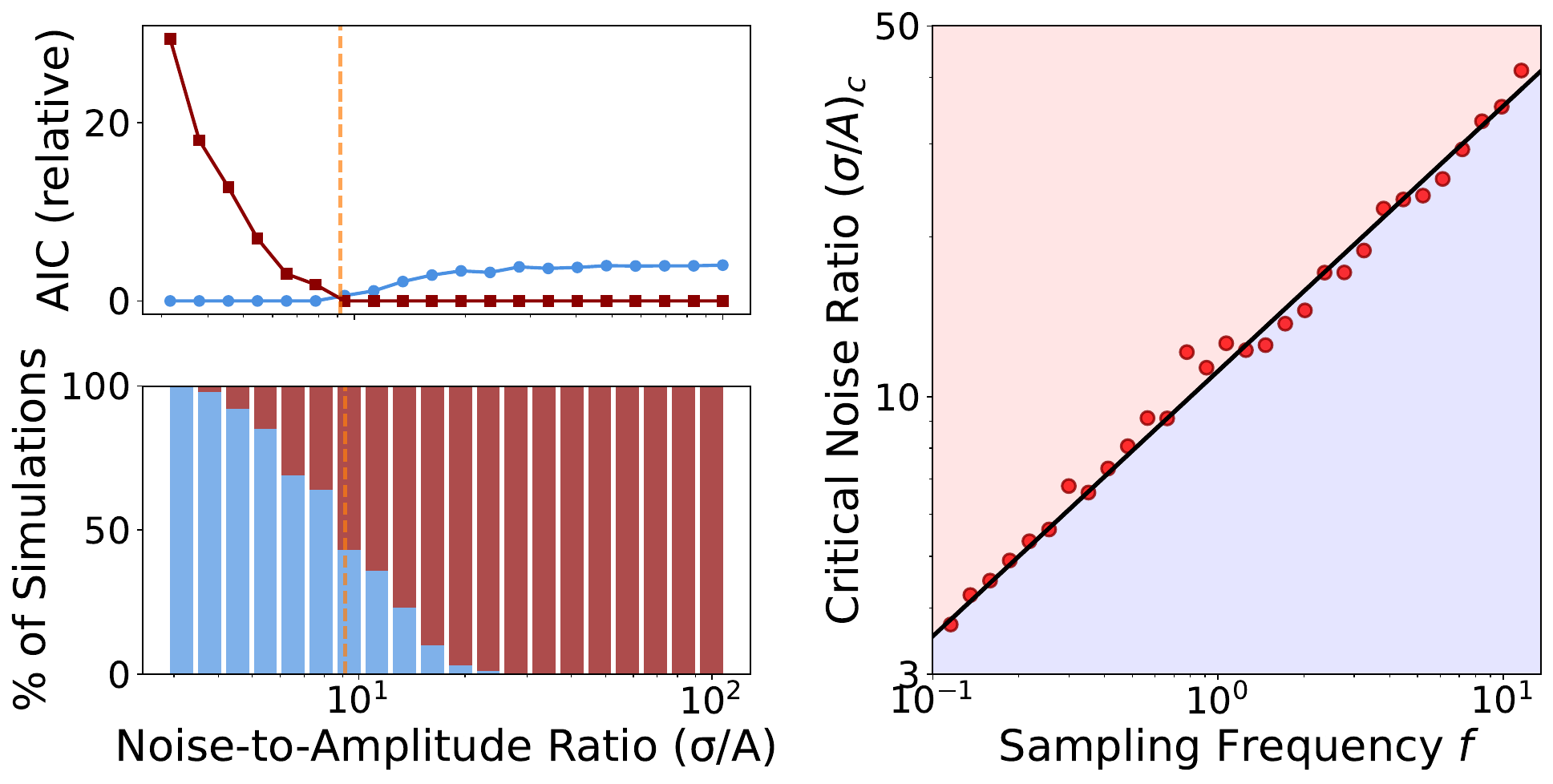}
    \caption{\textbf{Oscillation: Selected model varies with sampling rate.} 
    \textbf{Left panels:} Model selection at a single sampling frequency ($f = 0.67$) showing (top) normalized AIC values for the SHO model (blue circles) and pure noise model (red squares) as functions of noise-to-amplitude ratio $\sigma/A$, and (bottom) percentage of Monte Carlo simulations preferring each model. Vertical dashed line indicates the theoretical crossover point $(\sigma/A)_c = \sqrt{f t_{\max}/8}$. 
    \textbf{Right panel:} Numerical confirmation of $\sqrt{f}$ scaling law critical noise ratio.  Red filled circles: simulation, black line: theory. Parameters: $(A, \omega, f_0, t_{\max}) = (1, 2\pi, 1, 1000)$. }
    \label{fig:sho_samplingfreq}
\end{figure}

As with the decay model, we examine noisy observations $x(t) \sim A \cos(\omega t) + \eta$, where $\eta \sim \mathcal{N}(0,\sigma^2)$, and ask whether information criteria correctly identify the true oscillatory model versus a null hypothesis of pure noise. Ideally, model selection should depend only on the signal-to-noise ratio and not on the experimenter's choice of sampling frequency. However, as we demonstrate, sampling rate fundamentally affects this comparison. For fixed number of samples $M$, oscillatory systems exhibit the same qualitative behavior as equilibrating systems: both undersampling and oversampling degrade model selection, with an optimal intermediate frequency.

For oscillatory systems with \textit{fixed observation window} $t_{\max}$, increasing the sampling frequency $f$ increases the number of samples $M = t_{\max} \cdot f$, providing more independent observations of the oscillation. In this regime, analysis shows that the critical signal-to-noise ratio for model selection scales as
\begin{equation}
\left(\frac{\sigma}{A}\right)_c \sim \sqrt{t_{\max} f}.
\label{eq:shonoisefreq}
\end{equation}
This $\sqrt{f}$ dependence means that higher sampling frequencies allow detection of oscillations even in noisier data, provided the number of samples $M$ is sufficient for stable parameter estimation ($M \gtrsim 10$). For the parameters in Fig.~\ref{fig:sho_samplingfreq}, with oscillation period 1 (frequency $f_0 = 1$) and $t_{\max} = 1000$, the critical ratios are approximately $(\sigma/A)_c \approx 3.5$ for $f/f_0 = 0.1$ and $(\sigma/A)_c \approx 35$ for $f/f_0 = 10$.

Importantly, this scaling holds when the number of samples $M \gg k$ is sufficient for stable parameter estimation. For small $M$ or fixed $M$ (where increasing $f$ decreases $t_{\max} = M/f$), oscillatory systems exhibit similar oversampling pathologies as equilibrating systems: sampling too finely over too short a window prevents detection of the oscillatory structure. Thus, the $\sqrt{f}$ scaling represents an additional regime available to oscillatory systems when the observation window is held constant, rather than a fundamental difference from equilibrating systems. In the regime where this scaling holds, oversampling provides substantially wider noise tolerance for correct model identification compared to undersampling, which makes oscillation detection difficult except at very high signal-to-noise ratio (SNR).

\textbf{Deterministic chaos:} A third major dynamical regime beyond equilibration and oscillation is deterministic chaos, exemplified by the Lorenz system \cite{lorenz1963mechanics, strogatz2000}. We analyze a single trajectory component of the Lorenz system $X(t)$ (using standard parameters $(\sigma, \rho, \beta) = (10, 28, 8/3)$, yielding Lyapunov exponent $\lambda \approx 0.906$) and take $\lambda^{-1}$ as the characteristic time scale and $A \approx 38$ (the rough attractor amplitude) as the characteristic amplitude. Noisy observations are modeled as $x(t) \sim X(t) + \eta$, $\eta \sim \mathcal{N}(0, \sigma^2)$, and we assess how often this ``true'' model is selected when compared with a pure noise null.

For fixed observation time $t_{\max}$, increasing sampling frequency $f$ raises the sample count $M = t_{\max} f$, enhancing statistical power. The critical noise level for chaos model selection appears to scale with frequency as $(\sigma/A)_c \propto f^{\alpha}$, $\alpha \approx 0.51$, very close to $\sqrt{f}$ as in the SHO example above. Thus, higher sampling improves noise tolerance without a high-frequency failure mode.

However, chaotic systems are slightly more sensitive to undersampling due to their broadband, multi-scale structure. When $f \ll \lambda$, chaos becomes nearly indistinguishable from noise unless observations are almost noise-free. The Lyapunov exponent $\lambda$ therefore serves as a practical lower bound for adequate sampling, though the exact relationship among $f$, noise tolerance, and model selection remains empirical.

\section{Dimension Dependence}
\label{sec:dimdep}

We revisit the equilibrating system Eq.~\eqref{eq:decay} in a higher dimensional context, considering an $N$-dimensional decay process $x_{i}(t) = x_i(0)e^{-\lambda t} + \eta_{i}$, $\eta_{i}\sim\mathcal{N}(\mu,\sigma^{2})$, sampled as before at fixed frequency $f$.\footnote{Real-world systems might be expected to have dynamics coupled across dimensions, however we choose this example as a simplest case that may also be taken to represent the slowest-mode relaxation to equilibrium along an eigendirection.}

Our analysis depends critically on how data collection scales with system size $N$ and model complexity. We consider three plausible scenarios: \textbf{(1) Fixed data per dimension:} $M$ constant, regular sampling; 
\textbf{(2) Fixed total effort:} $M \cdot N$ constant, e.g., limited total storage available;
\textbf{(3) Combinatorial effects:} $M/N$ constant, i.e., total data $\sim N^2$ due to, e.g., network effects in dynamics.

\begin{figure}[t]
  \centering
  \includegraphics[width=\linewidth]{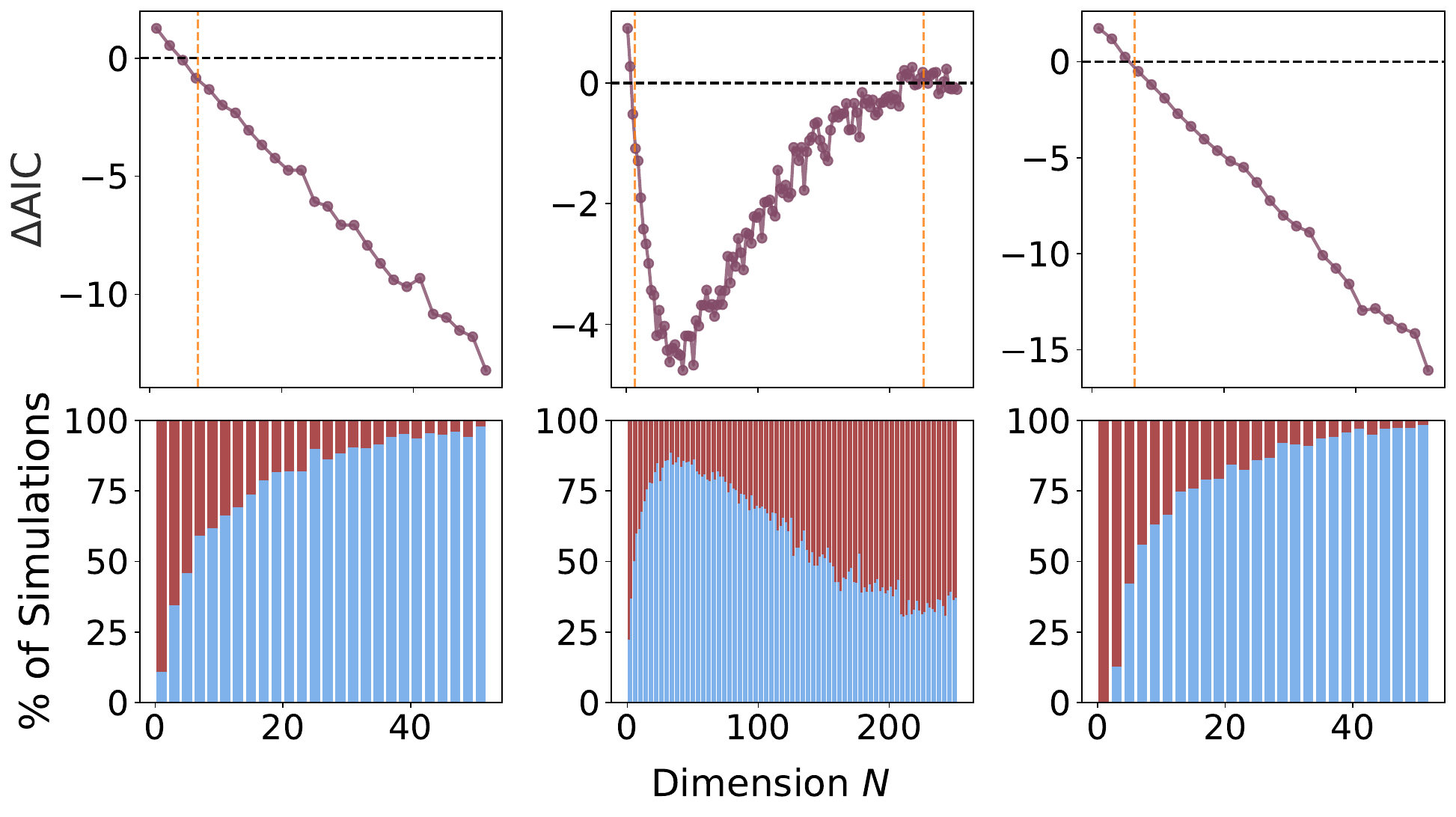}
  \caption{\textbf{Equilibration: Selected model varies with system dimension.} 
\textbf{Top:} $\Delta\text{AIC} = \text{AIC}_{\text{decay}} - \text{AIC}_{\text{noise}}$. Vertical orange dashed lines indicate theoretical crossover points $\Delta\text{AIC} = 0$. \textbf{Bottom:} Proportion of simulations (500 trials) selecting equilibration (blue) or pure noise (red) model. 
\textbf{Left panels:} Fixed data per dimension ($M = 100$ constant). Single crossover at $N_{\text{crit}} \approx 7$. 
\textbf{Middle panels:} Fixed total data ($MN = 1250$ constant). Two crossovers predicted: $N_{\text{crit}}^{\text{(low)}} \approx 6$ (large $M$ per agent) and $N_{\text{crit}}^{\text{(high)}} \approx 225$ (small $M$ per agent). 
\textbf{Right panels:} Combinatorial effects ($M/N = 10$ constant). Single crossover at $N_{\text{crit}} \approx 6$. 
All panels have fixed sampling frequency $f = 10$. Parameters: $x_0 = 2$ (known a priori, not fitted), ($\lambda$, $\mu$, $\sigma) = (1, 1, 8)$.}
\label{fig:dim_dependence_fixeddt}
\end{figure}

Fig.~\ref{fig:dim_dependence_fixeddt} demonstrates markedly different behaviors across scenarios. \textbf{Case 1:} Increasing $N$ provides more independent evidence while maintaining statistical power per agent; \textbf{Case 2:} Resource constraints create a trade-off such that the selected model may go from null to decay and back to null as dimension increases; \textbf{Case 3:} Data collection scaling with complexity amplifies the benefits of increased dimensionality.

To explain the crossovers in Fig.~\ref{fig:dim_dependence_fixeddt}, we compare the decay model with a pure noise null model analytically. We define variance 
\[
    S_i = M^{-1}\sum_{j=1}^{M}\left[x_i(0)e^{-\lambda t_{j}}-\musig_i \right]^{2},
\]
where $\musig_i = M^{-1}\sum_{j=1}^M x_i(0)e^{-\lambda t_j}$.  For simplicity, we first consider the case where all initial conditions are known (not fitted) and uniform and obtain
\begin{equation}\label{eq:dAIC_general}
  \Delta\mathrm{AIC}=2-NM\ln\!\Bigl(1+\frac{S}{\sigma^{2}}\Bigr),
\end{equation}
where $S_i = S\, \forall \, i$ is universal in this case (alternative assumptions about initial conditions are explored in SM).  Setting $\Delta\mathrm{AIC} = 0$ then yields: \begin{enumerate}
    \item ($M$ constant), $\Ncrit^{(M)} = 2/[M\ln(1 + S/\sigma^2)]$;
    \item ($MN = C$), $N = 2/[C\ln(1 + S/\sigma^2)]$ (here $S$ depends on $M = C/N$ so $\Ncrit^{(MN)}$ is implicit; two positive solutions exist);
    \item ($M/N = C$), $N = \sqrt{2/[C\ln(1 + S/\sigma^2)]}$ (here also $S$ depends on $M = C N$; only one positive solution for $\Ncrit^{(M/N)}$ exists).
\end{enumerate}  
See Appendix \ref{app:dimdep} for detailed calculations and asymptotic limits.

One takeaway from this examination of the impact of dimension on model selection is that, once again, the model selected may change in surprising ways.  Measuring one noisy relaxation process could lead the observer to believe it is best explained by pure noise, while measuring 100 could reverse that conclusion, while measuring 1000 could again lead to the pure noise conclusion---even when exactly the same dynamics govern every process.

\section{Discussion and Conclusions}

We have attempted to illustrate some weaknesses of model selection via information criteria in the context of dynamical systems.  Our crossover formulas make explicit predictions for particular dynamical systems, showing how the selected model may depend on sampling frequency and system dimension.

The fundamental issue is that AIC, BIC, and related criteria assume independent observations, an assumption violated by nearly any dynamical system. The systematic failures we highlight potentially affect thousands of published papers across physics. Our work provides an illustrative analytical characterization for first-order linear relaxation and harmonic oscillation---two ubiquitous dynamical motifs in physics. We expect that our results will help researchers compute effective sampling protocols given estimated system parameters, and encourage caution when modeling with dynamical data.

\bibliographystyle{apsrev4-2}
\bibliography{references}

\clearpage
\newpage

\setcounter{secnumdepth}{2}
\appendix

\section{Sampling Frequency Dependence: Detailed Derivations} 
\label{app:freqcrossover}

\subsection{Exponential Decay Model}
The decay model has parameters $\{x_0, \lambda, \mu, \sigma\}$ where $\mu$ is the noise mean and $\sigma$ is the noise standard deviation. The noise model has parameters $\{\mu, \sigma\}$. Assuming Gaussian noise,
\[
    \ln L = -\frac{M}{2}\ln(2\pi\sigma^2)-\frac{1}{2\sigma^2}\sum_{j=1}^M r_j^2,
\]
with residuals $r_j=x_j-\hat x_j$.  For the decay model, $r_j=x_j-\hat x_0e^{-\hat\lambda\,j/f} - \hat\mu$; for the noise model, $r_j=x_j-\hat\mu$. Hence
\begin{align}
    \mathrm{AIC}_{\textrm{decay}} &=8
    + \frac{1}{\hat\sigma^2}\sum_{j=1}^M\bigl(x_j-\hat x_0e^{-\hat\lambda\,j/f}-\hat\mu\bigr)^2 \nonumber \\
    & \qquad \qquad + M\ln(2\pi\hat\sigma^2),
    \label{aicdecay}\\
    \mathrm{AIC}_{\textrm{noise}}
    &=4
    +\frac{1}{\hat\sigma^2}\sum_{j=1}^M(x_j-\hat\mu)^2
    +M\ln(2\pi\hat\sigma^2).
    \label{aicnoise}
\end{align}
The AIC difference is then
\begin{equation}
    \Delta\mathrm{AIC} = 4 + 2(\mathrm{NLL}_{\textrm{decay}} - \mathrm{NLL}_{\textrm{noise}}),
\end{equation}
where NLL denotes negative log-likelihood. When both models fit the noise parameters well ($\hat\sigma_{\textrm{decay}} \approx \hat\sigma_{\textrm{noise}} \approx \sigma$), this simplifies to
\begin{equation}
\Delta\mathrm{AIC} \approx 4 - \frac{\Delta\mathrm{RSS}}{\sigma^2},
\end{equation}
where $\Delta\mathrm{RSS} = \mathrm{RSS}_{\textrm{noise}} - \mathrm{RSS}_{\textrm{decay}}$ is the difference in residual sum of squares.

For data generated from the true decay model with parameters $\{x_0, \lambda\}$ and noise $\eta \sim \mathcal{N}(0, \sigma^2)$, we can show that
\begin{equation}
    \Delta\mathrm{RSS} = M x_0^2 \mathrm{Var}[e^{-j\lambda/f}],
\end{equation}
where the variance is computed over the sample indices $j=0,1,...,M-1$.

\subsubsection{Low-frequency crossover (\texorpdfstring{$f \ll \lambda$}{f << lambda})}

In this regime, $\Lambda = \lambda M/f \gg 1$, meaning the signal decays completely over the observation window. Only the first $\sim f/\lambda$ points contain significant signal. The variance calculation yields
\begin{equation}
\mathrm{Var}[e^{-j\lambda/f}] \approx \frac{f}{2M\lambda},
\end{equation}
giving
\begin{equation}
\Delta\mathrm{RSS} \approx \frac{x_0^2 f}{2\lambda}.
\end{equation}

Setting $\Delta\mathrm{AIC} = 0$ for the crossover, we obtain
\begin{equation}
    f_c^{(1)} = \frac{8\lambda\sigma^2}{x_0^2}.
\end{equation}

Below this frequency, sparse sampling captures insufficient decay structure to justify the extra two parameters.

\subsubsection{High-frequency crossover (\texorpdfstring{$f \gg \lambda$}{f >> lambda})}

In this regime, $\Lambda = \lambda M/f \ll 1$, meaning the signal changes little over the observation window. Expanding $e^{-j\lambda/f} \approx 1 - j\lambda/f$, the variance of a linear trend gives
\begin{equation}
    \mathrm{Var}[e^{-j\lambda/f}] \approx \frac{\Lambda^2}{12} = \frac{\lambda^2 M^2}{12f^2},
\end{equation}
yielding
\begin{equation}
    \Delta\mathrm{RSS} \approx \frac{x_0^2 M^3 \lambda^2}{12f^2}.
\end{equation}

Setting $\Delta\mathrm{AIC} = 0$ for the crossover, we obtain
\begin{equation}
    f_c^{(2)} = \frac{M^{3/2}\lambda x_0}{4\sqrt{3}\sigma}.
\end{equation}

Above this frequency, the observation window is too short to capture sufficient decay, and the signal appears nearly constant.

\subsection{Simple Harmonic Oscillator}
For the SHO model we proceed similarly as in Appendix \ref{app:freqcrossover}. We have signal $s_j = A\cos(j\omega/f)$, which yields $\Delta\mathrm{RSS} = M \mathrm{Var}[A \cos(j\omega/f)]$.  Defining $\phi = \omega/f$ as the phase increment per sample, we compute the mean as
    $\frac{A}{M} \sin(M\phi/2) \cos\left[(M-1)\phi/2\right] / \sin(\phi/2)$ 
and second moment 
    $\frac{1}{2}A^2 + \frac{1}{2}A^2 M^{-1} \sin(M\phi) \cos[(M-1)\phi] / {\sin(\phi)}$.
When the number of samples is sufficiently large ($M \gg 4$ (number of parameters), and assuming $\phi$ is not commensurate with the oscillation frequency to avoid aliasing), both correction terms are $O(1/M)$, yielding ${\mathrm{Var}[A\cos(j\phi)] \approx A^2/2}$---approximately independent of both $f$ and $\omega$ for $M \gtrsim 10$. Using $M = t_{\max} f$, we find $\Delta\mathrm{RSS} \approx t_{\max} f A^2/2$ and $\Delta\mathrm{AIC} \approx 4 - t_{\max} f A^2/(2\sigma^2)$. This yields the crossover frequency ${f_c = 8\sigma^2/(t_{\max} A^2)}$, or, equivalently, ${(\sigma/A)_c = \sqrt{t_{\max} f/8}}$.

\section{Dimension Dependence: Detailed Derivations}
\label{app:dimdep}

\noindent \textbf{Signal Variance Function:}
For exponential decay $x(t) = x(0)e^{-\lambda t}$ sampled at $M$ points with spacing $\Delta t$, the signal variance is $S = x(0)^2 g(\alpha)$ where $\alpha = \lambda t_{\max} = \lambda M\Delta t$ is the dimensionless observation time and
\begin{equation}\label{eq:g_function}
    g(\alpha) = \frac{1-e^{-2\alpha}}{2\alpha} - \left(\frac{1-e^{-\alpha}}{\alpha}\right)^2.
\end{equation}
For short observation ($\alpha \ll 1$), Taylor expansion gives $g(\alpha) \approx \alpha^2/12$. For long observation ($\alpha \gg 1$), exponential terms vanish yielding $g(\alpha) \approx 1/(2\alpha)$. The detailed calculations can be found in SM.

\vspace{5mm}
\noindent \textbf{Case 1: Fixed \texorpdfstring{$M$, } (Fixed \texorpdfstring{$\Delta t$}):}
With $M$, $\Delta t$, and thus $\alpha = \lambda M\Delta t$ all constant, $S$ remains constant. From $\Delta\mathrm{AIC} = 2 - NM\ln(1 + S/\sigma^2) = 0$:
\begin{equation}
    N_{\textrm{crit}}^{(M)} = \frac{2}{M\ln\left(1 + \frac{S}{\sigma^2}\right)}, \quad S = x(0)^2 g(\lambda M\Delta t).
\end{equation}

\noindent \textbf{Case 2: Fixed \texorpdfstring{$MN = C$, } (Fixed \texorpdfstring{$\Delta t$})}

Here $M = C/N$ decreases with $N$, making $t_{\max} = C\Delta t/N$ and $\alpha(N) = \lambda C\Delta t/N$ both $N$-dependent. Thus $S(N) = x(0)^2 g(\lambda C\Delta t/N)$ varies, yielding:
\begin{equation}
\Delta\mathrm{AIC} = 2 - C\ln\left[1 + \frac{x(0)^2 g(\lambda C\Delta t/N)}{\sigma^2}\right].
\end{equation}

\textbf{Low $N$ (large $M$, large $\alpha$):} Using $g(\alpha) \approx 1/(2\alpha)$ for large $\alpha$:
\begin{equation}
N_{\textrm{crit}}^{({\textrm{low}})} \approx \frac{2}{C\ln\left(1 + \frac{x(0)^2}{2\lambda C\Delta t\sigma^2}\right)}.
\end{equation}

\textbf{High $N$ (small $M$, small $\alpha$):} Using $g(\alpha) \approx \alpha^2/12$ and $\ln(1+x) \approx x$ for small argument:
\begin{equation}
\Delta\mathrm{AIC} \approx 2 - \frac{Cx(0)^2(\lambda C\Delta t)^2}{12\sigma^2 N^2}.
\end{equation}
Setting to zero:
\begin{equation}
N_{\textrm{crit}}^{({\textrm{high}})} = \frac{\lambda C\Delta t \cdot x(0)}{2\sqrt{6}\sigma}\sqrt{C}.
\end{equation}

Case 2 exhibits two crossovers: decay model preferred in intermediate range $N_{\textrm{crit}}^{({\textrm{low}})} < N < N_{\textrm{crit}}^{({\textrm{high}})}$.

\vspace{5mm}
\noindent \textbf{Case 3: Fixed \texorpdfstring{$M/N = C$, } (Fixed \texorpdfstring{$\Delta t$})}

With $M = CN$ increasing linearly with $N$, we have $t_{\max} = CN\Delta t$ and $\alpha(N) = \lambda CN\Delta t$ both growing with $N$:
\begin{equation}
\Delta\mathrm{AIC} = 2 - CN^2\ln\left(1 + \frac{x(0)^2 g(\lambda CN\Delta t)}{\sigma^2}\right).
\end{equation}

For large $N$, using $g(\lambda CN\Delta t) \approx 1/(2\lambda CN\Delta t)$:
\begin{equation}
N_{\textrm{crit}}^{(\textrm{M}/\textrm{N})} = \sqrt{\frac{2}{C\ln\left(1 + \frac{x(0)^2}{2\lambda CN\Delta t\sigma^2}\right)}}.
\end{equation}
This is implicit since $S$ depends weakly on $N$ for large $N$, but provides good approximation.


\newpage \clearpage


\section*{Supplementary material}

# Supplemental Material: Information Criteria Fail for Dynamical Systems

## I. ADDITIONAL ANALYSIS: SAMPLING FREQUENCY DEPENDENCE

### A. Equilibrating Systems

**Fixed $f$ case:** When sampling frequency is fixed, increasing the number of samples improves model selection, as shown in Fig. S1. This happens because, as $M$ increases, so does $t_{\max}$. Hence, we are able to recover the true model when $t_{\max}$ is long enough to capture sufficient decay structure. The figure also confirms the validity of our theoretical crossover dependence on $M$, as computed in Eq. 4 of the main text.

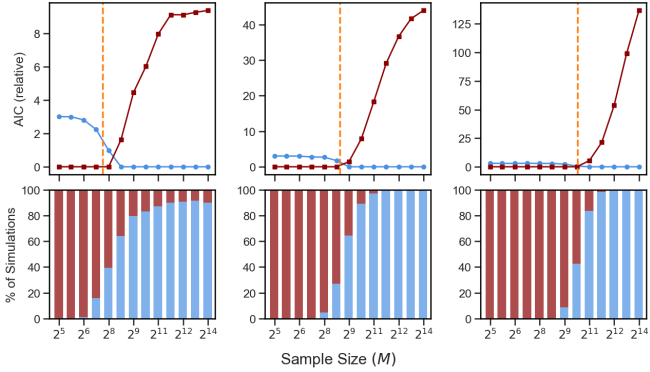

Figure S1. **Equilibration: Selected model varies with sample size.** Comparison of AIC values for both true model (blue circles) and pure noise model (red squares) as a function of sample size $M$. Left, middle, and right upper panels correspond to sampling frequencies $f = 2^4, 2^6$, and $2^8$, respectively. Lower panels show the proportion of simulations in which each model is selected (1000 trials each). The theoretical crossover (orange) occurs when the observation window becomes long enough to capture decay dynamics. From the high-frequency crossover condition, the critical sample size scales as $M_c \sim (f\sigma/\lambda x_0^2)^{2/3}$. Parameters: $(x_0, \lambda, \sigma) = (1, 0.1, 2.5)$.

**Summary:** In Table S1 we summarize all our findings from the main text and the above subsection regarding model selection in equilibrating systems and the dependence on sampling frequency.

### B. Oscillatory Systems

**Distribution mismatch caveat:** The likelihood analysis discussed in the main text assumes a Gaussian distribution for the observed data. However, for low noise, the distribution of observed values approaches $P(x) \propto 1/\sqrt{A^2 - x^2}$, very distinct from a Gaussian. This distribution mismatch affects not only the model selec-

| Regime | Condition | AIC Behavior |
|---|---|---|
| Low freq. | $f \ll f_c^{(1)}$ | Insufficient time resolution. Increasing $M$ doesn't help. Noise model preferred. |
| Mid freq., low $M$ | $f_c^{(1)} \lesssim f \lesssim f_c^{(2)}$, $M \ll M_c$ | Sampling rate is sufficient, but the total duration $t_{\max} = M/f$ is too short to capture decay. Noise model preferred. |
| Mid freq., high $M$ | $f_c^{(1)} \lesssim f \lesssim f_c^{(2)}$, $M \gtrsim M_c$ | Sampling and observation window are sufficient to resolve decay. Decay model preferred. |
| High freq. | $f \gg f_c^{(2)}$ | Duration $t_{\max} = M/f \ll \lambda^{-1}$. Samples are too closely spaced in time. Noise model preferred. |

Table S1. Model selection behavior vs. sampling frequency $f$, number of samples $M$, and system parameters. Here, $f_c^{(1)} = 8\lambda\sigma^2/x_0^2$, $f_c^{(2)} = (M^{3/2}\lambda x_0^2)/(4\sqrt{3}\sigma)$, and $M_c = (4\sqrt{3}f\sigma/\lambda x_0^2)^{2/3}$.

tion but also parameter estimation quality and can shift effective crossover points. This effect is even more pronounced as noise grows. We show visually how the distribution changes as noise level grows in Fig. S2.

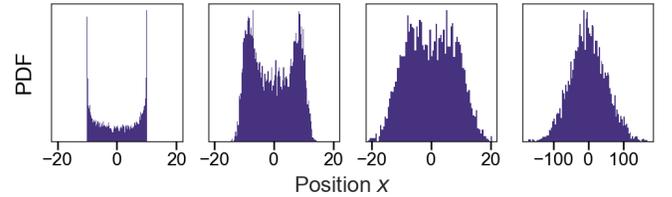

Figure S2. **Oscillation: distribution of data vs noise level.** Changing data distribution for the simple harmonic oscillator $x(t) \sim 10\cos(t) + \eta$ with varying $\sigma$ ($\sigma = 0.01, 1.5, 4, 50$ from left to right).

### C. Deterministic Chaos

We investigated deterministic chaos in the Lorenz system primarily through computational experiment, since theory for model selection in this case is challenging. However, it appears that the Lorenz system shares much in common, at least qualitatively, with oscillatory systems, as shown in Figs. S3 and S4. As with the SHO, increasing sampling frequency $f$ (for fixed observation

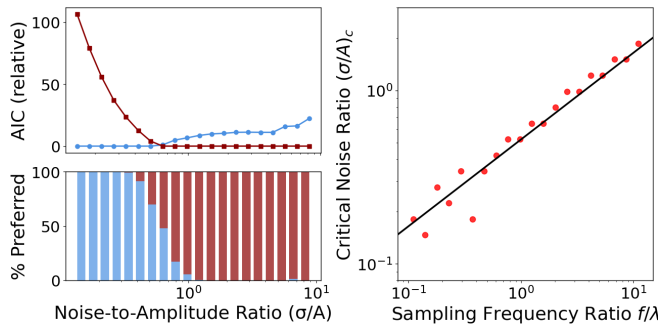

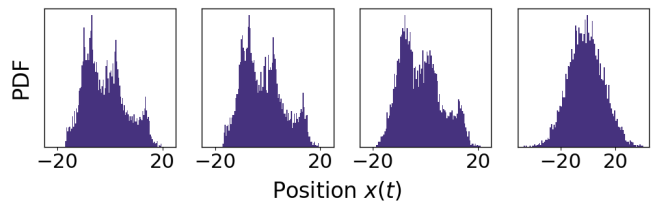

Figure S3. **Chaos: Model selection depends on sampling rate and noise. Left panels:** Model selection at a single sampling frequency showing (top) relative AIC values for the Lorenz model (blue circles) and white noise model (red squares) as functions of noise-to-amplitude ratio $\sigma/A$, and (bottom) percentage of Monte Carlo simulations preferring each model. The crossover occurs near $\sigma/A \sim 1$ where models become equally preferred. **Right panel:** Empirically determined critical noise levels where models are equally preferred (50% crossover point) as a function of sampling frequency. Red circles show empirical measurements across a range of sampling frequencies. Reference line is $(\sigma/A) \sim \sqrt{t_{\max} f}$, the analytical insight drawn from the simple harmonic oscillator (Eq. 6 example of the main text). Parameters: $(\sigma, \rho, \beta) = (10, 28, 8/3)$, $\lambda \approx 0.906$, $A \approx 38$, $t_{\max} = 100$. For the left panel: $f = 1$ ($f/\lambda \approx 1.1$), giving $M = 100$ samples.

time $t_{\max}$) increases the number of samples $M = t_{\max} f$, providing more statistical power to distinguish chaotic structure from noise. The critical noise level for chaos detection increases with sampling frequency, meaning that higher sampling rates allow identification of chaotic dynamics even in noisier data.

The right panel of Fig. S3 shows the empirical relationship between sampling frequency and critical noise level. The Lorenz system displays a very similar dependence on sampling frequency as SHO, following an apparent power law that appears consistent with $(\sigma/A)_c \propto \sqrt{f}$ as was shown for the SHO. This scaling suggests that a chaotic signal, at least in the Lorenz system example, behaves similarly (with regard to sampling rate dependence in model selection) to an oscillatory signal (though we do perhaps see notable deviations for very low sampling frequencies).

The physical origin of this sensitivity to noise across different sampling regimes likely arises from the multi-scale temporal structure inherent in chaotic dynamics. In systems such as the Lorenz attractor, exponential divergence of nearby trajectories produces correlations that extend across sampling intervals in ways fundamentally different from those in simple periodic oscillations. The resulting distribution mismatch is shown in Fig. S4, and the discussion in the previous section applies here as well. Whereas a sinusoidal signal concentrates its power at a single frequency, chaotic attractors exhibit broad spectral content, and adequate temporal resolution is required to capture this richer dynamical structure.

Figure S4. **Chaos: distribution of data vs noise level.** Changing data distribution for the $x(t)$ trajectory of the Lorenz system with observation noise of varying standard deviation $\sigma$ ($\sigma = 0.01, 1.5, 4, 50$ from left to right).

## II. DETAILED DERIVATIONS FOR DIMENSION DEPENDENCE

### A. Signal Variance Function $g(\alpha)$

For an exponential decay signal $x(t) = x(0)e^{-\lambda t}$ sampled at $M$ equally-spaced time points $t_j = j\Delta t$ for $j = 0, 1, \ldots, M-1$, the signal variance is:

$$S = \frac{1}{M} \sum_{j=0}^{M-1} \left( x(0)e^{-\lambda t_j} - \mu_{\text{signal}} \right)^2,$$

where $\mu_{\text{signal}} = \frac{1}{M} \sum_{j=0}^{M-1} x(0)e^{-\lambda t_j}$ is the temporal mean. The following calculations find analytical expressions for the variance in specific asymptotic limits.

For large $M$, the sums can be approximated by integrals over the observation window $[0, t_{\max}]$ where $t_{\max} = (M-1)\Delta t \approx M\Delta t$:

$$\mu_{\text{signal}} \approx \frac{x(0)}{t_{\max}} \int_0^{t_{\max}} e^{-\lambda t}\, dt = \frac{x(0)}{t_{\max}} \frac{\left(1 - e^{-\lambda t_{\max}}\right)}{\lambda} \quad \text{(S1)}$$

$$S \approx \frac{x(0)^2}{t_{\max}} \int_0^{t_{\max}} e^{-2\lambda t}\, dt - \mu_{\text{signal}}^2$$

$$= \frac{x(0)^2}{t_{\max}} \frac{\left(1 - e^{-2\lambda t_{\max}}\right)}{2\lambda} - \left[ \frac{x(0)}{t_{\max}} \frac{\left(1 - e^{-\lambda t_{\max}}\right)}{\lambda} \right]^2. \quad \text{(S2)}$$

Defining the dimensionless observation time $\alpha = \lambda t_{\max}$, we can factor out $x(0)^2$ and write

$$S(\alpha) = x(0)^2\, g(\alpha),$$

where

$$g(\alpha) = \frac{1 - e^{-2\alpha}}{2\alpha} - \left( \frac{1 - e^{-\alpha}}{\alpha} \right)^2. \quad \text{(S3)}$$

**Short observation time ($\alpha \ll 1$):** Using Taylor expansions for $e^{-x}$ yields

$$\frac{1 - e^{-2\alpha}}{2\alpha} = 1 - \alpha + \frac{2\alpha^2}{3} + O(\alpha^3)$$

and

$$\left(\frac{1-e^{-\alpha}}{\alpha}\right)^2 = 1 - \alpha + \frac{7\alpha^2}{12} + O(\alpha^3).$$

Therefore

$$g(\alpha) \approx \frac{2\alpha^2}{3} - \frac{7\alpha^2}{12} = \frac{\alpha^2}{12} \quad \text{for } \alpha \ll 1. \quad \text{(S4)}$$

**Long observation time ($\alpha \gg 1$):** For large $\alpha$, exponential terms vanish:

$$g(\alpha) \approx \frac{1}{2\alpha} - \left(\frac{1}{\alpha}\right)^2 = \frac{1}{2\alpha} - \frac{1}{\alpha^2}.$$

The leading term dominates:

$$\boxed{g(\alpha) \approx \frac{1}{2\alpha} \quad \text{for } \alpha \gg 1.} \quad \text{(S5)}$$

### B. Alternative Initial Condition Assumptions

As mentioned in our main discussion, initial conditions require careful interpretation. In our calculations in the main text, we treated initial values of all trajectories as known. The other two possibilities we lay out here are: (1) Treating all initial conditions as unknown but universal (i.e., a single fitted parameter); or (2) Treating the initial condition for each dimension as a distinct fitted parameter (i.e., $N$ fitted parameters).

Treating all trajectories as sharing a single initial value as known a priori keeps the complexity penalty fixed, so added dimensions contribute coherently to evidence accumulation. The case discussed in the paper shows clear threshold behavior: when $M \ln(1 + S/\sigma^2)$ is large enough, increasing $N$ strengthens recovery of the true model; otherwise the model is already disfavored at $N = 1$ and cannot be recovered at any dimension. If the initial conditions are universal and fitted, an extra penalty is incurred and we need to go to an even higher dimension to improve identifiability.

By contrast, if each trajectory has its own independently fitted initial condition, the number of free parameters grows with $N$, reversing the trend: as system size increases, the comparison shifts systematically toward the simpler noise model. This sharp contrast—from unconditional improvement, to conditional thresholds, to complete reversal—shows that the conventions used for initial-condition parameterization entirely determine whether collective evidence aids or hinders model identification.

**Independent Initial Conditions Analysis:** When each agent's initial condition is treated as an independent fitted parameter, the number of system parameters for the true model grows linearly as the dimension grows.

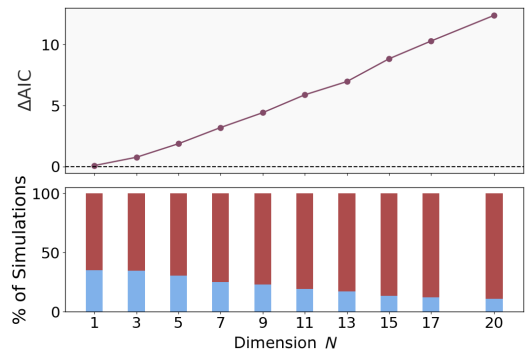

Figure S5. **Curse of dimensionality when initial conditions are independent parameters.** Top: $\Delta\text{AIC} = \text{AIC}_{\text{decay}} - \text{AIC}_{\text{noise}}$. Bottom: Proportion of simulations (500 trials) selecting equilibration (blue) or pure noise (red) model. Unlike the analysis related to Fig. 3, where increasing dimensionality enabled true model recovery at large $N$, here increasing dimension makes true model recovery less likely. The noise model is selected with increasing frequency as $N$ grows, reaching nearly 100% selection at higher dimensions. Parameters: $M = 1000$ fixed data points per dimension, $f = 10$, $x_i(0) \sim \mathcal{N}(0.7, 0.1^2)$, $\lambda = 1$, $\sigma^2 = 3$.

This leads the penalty for the true model to also grow linearly with dimension (see Fig. S5).

For equilibrating systems with fixed $M$, the calculations are as follows. We consider the means of the noise to be zero, for simplicity. The symbols are the same as in the main text section discussion.

**Parameter counts:**

- Decay model: $x_1(0), x_2(0), \ldots, x_N(0), \lambda, \mu, \sigma \to k_{\text{decay}} = N + 3$

- Noise model: $\mu, \sigma \to k_{\text{noise}} = 2$

- $\Delta k = N + 1$

**AIC difference:**

$$\Delta\text{AIC} = 2(N+1) - NM \ln\left(1 + \frac{S}{\sigma^2}\right)$$

$$= 2N + 2 - NM \ln\left(1 + \frac{S}{\sigma^2}\right)$$

$$= N\left[2 - M \ln\left(1 + \frac{S}{\sigma^2}\right)\right] + 2. \quad \text{(S6)}$$

**Slope analysis:** The slope with respect to $N$ is:

$$\frac{d(\Delta\text{AIC})}{dN} = 2 - M \ln\left(1 + \frac{S}{\sigma^2}\right).$$

- If $M \ln(1 + S/\sigma^2) > 2$: Negative slope → dimensionality blessing

- If $M \ln(1 + S/\sigma^2) < 2$: Positive slope → dimensionality curse

- At $N = 1$: $\Delta\text{AIC}(1) = 4 - M\ln(1 + S/\sigma^2)$

**Key observation:** If the slope were positive ($M\ln(1+S/\sigma^2) < 2$), then automatically $\Delta\text{AIC}(1) = 4 - M\ln(1+S/\sigma^2) > 2 > 0$, meaning the noise model would already be preferred at $N = 1$ and increasingly so for larger $N$. There is no crossover possible. Note, however, that for other cases where $M$ varies with $N$, the analysis is more complicated and crossovers may be possible.

### C. Fixed $t_\text{max}$ Cases

In our analysis in the main text, we introduced three cases (Case 1: Fixed $M$, Case 2: Fixed $MN$, Case 3: Fixed $M/N$). We made the additional assumption of fixed sampling frequency $f$ (equivalently, fixed $\Delta t$). Here we discuss the same cases, but with the assumption that $t_\text{max}$ is fixed instead of $f$. When observation time $t_\text{max}$ is held constant, $\alpha = \lambda/f$ is also constant, making $S = x(0)^2 g(\alpha)$ independent of $N$ (see Appendix B). Results for this case are visualized in Fig. S6 and are detailed below.

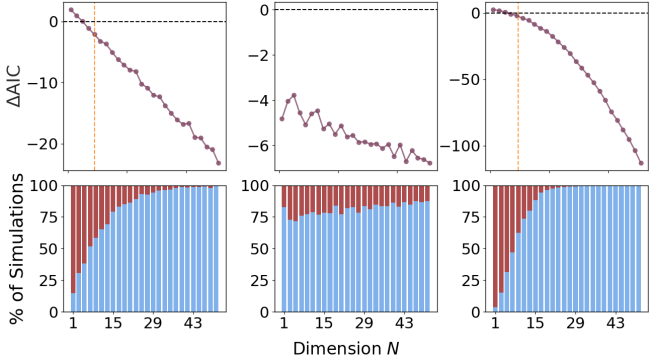

Figure S6. **Equilibration: Selected model varies with system dimension, fixed $t_\text{max}$. Top:** $\Delta\text{AIC} = \text{AIC}_\text{decay} - \text{AIC}_\text{noise}$. Vertical orange dashed lines indicate theoretical crossover points $\Delta\text{AIC} = 0$. **Bottom:** Proportion of simulations (500 trials) selecting equilibration (blue) or pure noise (red) model. **Left panels:** Fixed data per dimension ($M = 100$ constant). Single crossover at $N_\text{crit} \approx 9$. **Middle panels:** Fixed total data ($MN = 1500$ constant). No crossover predicted. **Right panels:** Combinatorial effects ($M/N = 10$ constant). Single crossover at $N_\text{crit} \approx 10$. All panels have fixed sampling frequency $f = 10$. Parameters: $x_0 = 2$ (known a priori, not fitted), $(\lambda, \mu, \sigma) = (1, 1, 6)$.

**Case 4: Fixed $M$, Fixed $t_\text{max}$.** Since $t_\text{max} = M/f$ is constant and $M$ is constant, $f$ must also be constant here. This is identical to Case 1 from the main text, in the fixed $f$ discussion.

**Case 5: Fixed $MN = C$, Fixed $t_\text{max}$.** With $M = C/N$ and $t_\text{max} = M/f$ constant:

$$f = \frac{M}{t_\text{max}} = \frac{C}{Nt_\text{max}}.$$

Thus $f$ varies with $N$, but $S$ remains constant (since $\alpha = \lambda t_\text{max}$ is constant), and

$$\Delta\text{AIC} = 2 - MN\ln\left(1 + \frac{S}{\sigma^2}\right) = 2 - C\ln\left(1 + \frac{S}{\sigma^2}\right).$$

This is **independent of** $N$, and there is no crossover with respect to dimension. A crossover does, however, occur with respect to $C$:

$$C_\text{crit} = \frac{2}{\ln\left(1 + \frac{S}{\sigma^2}\right)}, \quad S = x(0)^2 g(\lambda t_\text{max}). \quad \text{(S7)}$$

**Case 6: Fixed $M/N = C$, Fixed $t_\text{max}$:** With $M = CN$ and $t_\text{max} = M/f$ constant:

$$f = \frac{CN}{t_\text{max}}.$$

As above, $S$ is constant, so

$$\Delta\text{AIC} = 2 - CN^2 \ln\left(1 + \frac{S}{\sigma^2}\right).$$

Setting this to zero, we find

$$N_\text{crit} = \sqrt{\frac{2}{C\ln\left(1 + \frac{S}{\sigma^2}\right)}}. \quad \text{(S8)}$$

### D. Summary Table

We summarize our findings for dimension-dependence in all cases discussed here and in the main text in Table S2.

| Case | $S$ behavior | $N_\text{crit}$ formula |
|---|---|---|
| 1: Fixed $M$, $f$ | Constant | $\frac{2}{M\ln(1+S/\sigma^2)}$ |
| 2: Fixed $MN$, $f$ | $\propto g\left(\frac{\lambda C}{Nf}\right)$ | Low: $\frac{2}{C\ln(1+S/\sigma^2)}$ <br> High: $\sqrt{\frac{Cx(0)^2(\lambda C)^2}{24(\sigma f)^2}}$ |
| 3: Fixed $M/N$, $f$ | $\propto g\left(\frac{\lambda CN}{f}\right)$ | $\sqrt{\frac{2}{C\ln(1+S/\sigma^2)}}$ |
| 4: Fixed $M$, $t_\text{max}$ | Constant | $\frac{2}{M\ln(1+S/\sigma^2)}$ |
| 5: Fixed $MN$, $t_\text{max}$ | Constant | None |
| 6: Fixed $M/N$, $t_\text{max}$ | Constant | $\sqrt{\frac{2}{C\ln(1+S/\sigma^2)}}$ |

Table S2. Summary of critical dimension formulas for all cases with difference in number of parameters, $\Delta k = k_\text{decay} - k_\text{noise} = 1$. Note that formulas for $N_\text{crit}$ are exact when $S$ is constant, but approximate when not.

## III. BIC AND OTHER INFORMATION CRITERIA

The pathologies demonstrated for AIC extend directly to BIC, as both the methods share the same likelihood-based structure but with different complexity penalties.

For BIC, the penalty term scales as $k \ln N$ rather than $2k$, where $N$ is the number of data points. This similarity ensures they exhibit the same fundamental sampling rate and dimension dependence issues.

**Equilibrating Systems: Crossover Analysis.** For the equilibrating system example, the BIC difference is

$$\Delta \text{BIC} = \ln M - M \ln\left(1 + \frac{S}{\sigma^2}\right),$$

where we have used the fact that BIC has penalty $k \ln M$ and $\Delta k = 1$.

The BIC crossover frequencies can be derived by setting $\Delta \text{BIC} = 0$:

$$M_{\text{crossover}} = \exp\left(\frac{\ln M}{\ln(1 + S/\sigma^2)}\right).$$

This transcendental equation shows that BIC crossovers occur at different sampling frequencies than AIC crossovers, but the qualitative behavior remains the same: both under-sampling and over-sampling lead to incorrect model selection. We can see the same qualitative behavior as in our first section of the main paper in Fig. S7. We observe the same similarities between AIC and BIC behaviours across all experiments.

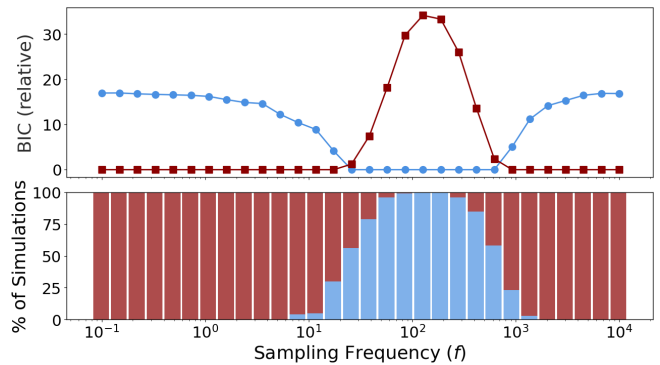

Figure S7. **Model Selection with BIC.** Top: Comparison of BIC values for true deterministic decay model (blue) and pure noise null model (red) as a function of sampling frequency $f$. Bottom: Proportion of trials in which each model is selected as a function of sampling frequency $f$ (100 trials). Note that the qualitative behaviour is exactly the same as our analysis with AIC. Parameters: $(x_0, \lambda, \mu, \sigma, M) = (1, 0.1, 0, 2.5, 5000)$. Points in upper panel represent means over trials at each frequency.

**Equilibrating Systems: Dimension Dependence Analysis.** For the multi-dimensional equilibrating system, we can derive the complete crossover analysis for BIC. Starting from the BIC difference:

$$\Delta \text{BIC} = \ln(MN) - NM \ln\left(1 + \frac{S}{\sigma^2}\right), \quad (S9)$$

Setting $\Delta \text{BIC} = 0$ for the crossover condition:

$$N_{\text{crit}}^{\text{BIC}} = \frac{\ln(MN)}{M \ln(1 + S/\sigma^2)}. \quad (S10)$$

This transcendental equation can be solved iteratively. For large $M$, we can approximate:

$$N_{\text{crit}}^{\text{BIC}} \approx \frac{\ln(M)}{M \ln(1 + S/\sigma^2)} + \frac{\ln(N_{\text{crit}}^{\text{BIC}})}{M \ln(1 + S/\sigma^2)}. \quad (S11)$$

The second term represents a logarithmic correction that becomes negligible for $M \gg 1$, yielding:

$$N_{\text{crit}}^{\text{BIC}} \approx \frac{\ln(M)}{M \ln(1 + S/\sigma^2)}. \quad (S12)$$

Comparing with the critical $N$ using AIC from dimension dependence section of the main text, we find:

$$\frac{N_{\text{crit}}^{\text{BIC}}}{N_{\text{crit}}^{\text{AIC}}} \approx \frac{\ln(M)}{2}.$$

Since $\ln(M) < 2$ for $M < e^2 \approx 7.4$, BIC requires fewer dimensions than AIC for small sample sizes, but more dimensions for large samples. This reflects BIC's stronger penalty for model complexity.



\end{document}